\documentclass[review, 3p, onecolumn]{elsarticle}
\usepackage{microtype}
\usepackage{graphicx}
\usepackage{amssymb}
\usepackage{amsmath}
\usepackage{lineno}
\usepackage{booktabs}
\usepackage{listings}
\usepackage{multirow}
\usepackage{lscape}
\usepackage{ bbold }
\usepackage{ esint }
\usepackage{algorithm}
\usepackage{algpseudocode}
\usepackage{enumitem}
\usepackage{siunitx}  

\usepackage{graphicx}    
\usepackage{subcaption}  

\PassOptionsToPackage{hyphens}{url}\usepackage[colorlinks=True]{hyperref}

\biboptions{sort&compress}

\usepackage{array}
\newcolumntype{L}[1]{>{\raggedright\let\newline\\\arraybackslash\hspace{0pt}}m{#1}}
\newcolumntype{C}[1]{>{\centering\let\newline\\\arraybackslash\hspace{0pt}}m{#1}}
\newcolumntype{R}[1]{>{\raggedleft\let\newline\\\arraybackslash\hspace{0pt}}m{#1}}

\def\equationautorefname~#1\null{(#1)\null}

\usepackage{framed} 
\usepackage{multicol} 
\usepackage{nomencl} 
\makenomenclature
\renewcommand*\nompreamble{\begin{multicols}{2}}
\renewcommand*\nompostamble{\end{multicols}}
\usepackage{etoolbox}
\renewcommand\nomgroup[1]{%
  \item[\bfseries
  \ifstrequal{#1}{L}{Latin letters}{%
  \ifstrequal{#1}{M}{Greek letters}{%
  \ifstrequal{#1}{F}{Calligraphic letters}{%
  \ifstrequal{#1}{S}{Subscripts and superscripts}{}}}}%
]}

\newcommand{\thermaldif}{\alpha}
\newcommand{\thermalcond}{\lambda}
\newcommand{\nbhe}{N_b}
\newcommand{\length}{L}
\newcommand{\positions}{\mathbf{p}}
\newcommand{\heatextr}{q}
\newcommand{\green}{\mathcal G}
\newcommand{\verticalkernel}{\mathcal Z}
\newcommand{\horizontalkernel}{\mathcal R}
\newcommand{\bheradius}{r_b}
\newcommand{\density}{\rho}
\newcommand{\flowrate}{\dot m}

\begin{document}

\begin{frontmatter}

\title{Optimization of Closed-Loop Shallow Geothermal Systems Using Analytical Models}

\author[MO]{Oliver Heinzel\corref{correspondingauthor}}
\ead{o.heinzel@tum.de}
\cortext[correspondingauthor]{Corresponding author}
\author[ENS]{Smajil Halilovic}
\author[ENS]{Thomas Hamacher}
\author[MO]{Michael Ulbrich}

\address[MO]{Technical University of Munich, Chair of Mathematical Optimization, Germany}
\address[ENS]{Technical University of Munich, Chair of Renewable and Sustainable Energy Systems, Germany}

\begin{abstract}
 Closed-loop shallow geothermal systems are one of the key technologies for decarbonizing the residential heating and cooling sector. The primary type of these systems involves vertical borehole heat exchangers (BHEs). During the planning phase, it is essential to find the optimal design for these systems, including the depth and spatial arrangement of the BHEs. In this work, we have developed a novel approach to find the optimal design of BHE fields, taking into account constraints such as temperature limits of the heat carrier fluid. These limits correspond to the regulatory practices applied during the planning phase. The approach uses a finite line source model to simulate temperature changes in the ground in combination with an analytical model of heat transport within the boreholes. Our approach is demonstrated using realistic scenarios and is expected to improve current practice in the planning and design of BHE systems.  
\end{abstract}


\end{frontmatter}

\begin{table*}[th!]
\footnotesize
\label{tab:nomenclature}
  \begin{framed}
  \nomenclature[M]{$\thermaldif$}{Thermal diffusivity [m$^2$/s]}
  \nomenclature[M]{$\thermalcond$}{Thermal Conductivity [W/(m K)]}
    \nomenclature[M]{$\beta$}{Normalized borehole conductance [-]}
       \nomenclature[M]{$\Psi$}{Local/Global Field Interaction term [K]}
        \nomenclature[M]{$\gamma$}{Borehole Characteristic parameter [-]}
  \nomenclature[L, B]{$\nbhe$}{Number of Boreholes [-]}
  \nomenclature[L, B]{R}{Resistance}
  \nomenclature[L, m]{$\flowrate$}{Flow rate [kg/s]}
    \nomenclature[L, T]{$T$}{Temperature [K]}
  \nomenclature[L, c]{$c$}{Specific heat capacity [J/(kg $\cdot$ K)]}
  \nomenclature[L, L]{$\length$}{Length of the Boreholes [m]}
  \nomenclature[L, x]{$x,y$}{Horizontal coordinates [-]}
  \nomenclature[L, z]{$z$}{Vertical coordinate [-]}
  \nomenclature[L, p]{$\positions$}{Positions of BHE [-]}
  \nomenclature[S]{$k,\hat k$}{Counter for BHEs}
  \nomenclature[S]{$r$}{Refrigerant}
  \nomenclature[L, q]{$\heatextr$}{Heat extraction rate per unit length [W/m]}
  \nomenclature[F]{$\green$}{Green's Function}
  \nomenclature[F]{$\verticalkernel$}{Vertical Kernel}
  \nomenclature[F]{$\horizontalkernel$}{Horizontal Kernel}
    \nomenclature[F]{$\mathcal{T}$}{Simulation end time [s]}
  \nomenclature[L, rb]{$\bheradius$}{Borehole radius}
  \nomenclature[M]{$\density$}{Density [kg/m$^3$]}
  \nomenclature[S]{inter}{pipe-to-pipe connection}
    \nomenclature[S]{s}{soil-to-pipe connection}
\nomenclature[S]{$^\circ$}{interior of a set}
    \printnomenclature
  \end{framed}
\end{table*}


\section{Introduction}


The growing global need for sustainable energy sources has significantly increased demand for closed-loop shallow geothermal systems over the last decade, particularly those with vertical borehole heat exchangers (BHEs) \cite{Halilovic.2022, Halilovic.2023}. These systems provide an excellent and environmentally friendly way to heat and cool buildings \cite{Rybach_2010, Spitler_2015}. However, their widespread use and economic viability often depend on keeping drilling costs low, which are directly related to the total BHE length, while simultaneously making sure the system is safe and works well throughout its lifespan \cite{Kavanaugh_2008, Hellstrom_1995}.

Planing an optimal BHE field is a challenging task which involves different areas of physics. Full-scale numerical models, like those that use finite element or finite difference methods \cite{Florides_2007, Halilovic.2022a}, can accurately simulate ground heat transfer, yet they are often too expensive to use in iterative optimization loops \cite{Eskilson_1987, Halilovic.2023b, Halilovic.2024}. To address these computational challenges, modular simulation environments like the Modelica Buildings Library \cite{WetterZuoNouiduiPang2014} have implemented efficient base classes for calculating thermal response factors, allowing for the dynamic simulation of large-scale borehole fields without the overhead of full-scale discretized models.
While temperature response functions, or g-functions, remain the only feasible method for long-term whole-building energy simulations, recent studies using open-source tools such as pygfunction \cite{CIMMINO2014641} have shown that managing the computational time and memory requirements for larger fields, often exceeding 100 boreholes, continues to present significant practical challenges \cite{osti_1737649}.
Moreover, traditional analytical solutions for BHE sizing frequently utilize notable simplifications that often neglect the complex, dynamic thermal coupling between the circulating heat carrier fluid and the neighboring soil \cite{Yavuzturk_1999}. These simplifications can make it more difficult to establish accurate predictions about how the system will work over time and may even lead to designs that do not meet operational needs over the system's lifetime. A significant limitation of existing design tools for BHEs is their restricted capability to accommodate complex, non-convex site geometries \cite{sanner1997earth, Peere2022}. Conventional methodologies predominantly rely on rigid rectangular or regular grid layouts, which often result in suboptimal land utilization or necessitate labor-intensive manual adjustments when confronted with irregular property boundaries, internal obstacles, or building footprints. To overcome these geometric constraints, this study explores a more flexible spatial optimization framework based on Centroidal Voronoi Tessellations (CVTs). Lloyd’s algorithm \cite{Lloyd1982} remains the foundational iterative procedure for approximating such tessellations; while originally developed for signal quantization, it has proven highly effective in modern practical applications requiring efficient spatial distribution, ranging from the autonomous deployment of multi-robot swarms \cite{Santos2019} to the coverage optimization of wireless sensor networks \cite{Zhu2021}.
In light of these difficulties, this paper seeks to create a computationally efficient and geometrically adaptable optimization framework for BHE field design. The specific goals of this work are threefold:
\begin{itemize} 
    \item  To make a quick, analytical model that can predict how well the BHE system will work over a 20-year period, specifically for the 1U-tube BHE setup, which is widely used. This model will be the main part of our optimization loop.
   \item To present a robust heuristic approach, derived from an adapted Lloyd's algorithm, capable of autonomously identifying optimal BHE coordinates within any non-convex domain, thereby maximizing inter-borehole spacing and distance from domain boundaries.
   \item  To combine the above parts into an optimization problem that minimizes the uniform length of BHEs while strictly following the critical fluid temperature limits ($T_{\min}$ and $T_{\max}$) for the whole simulation period, making sure that the heat pump works safely and efficiently.
\end{itemize}

\section{Mathematical Modeling}

 \subsection{Soil Temperature Model }

We examine the heat transfer within the subsurface, modeled as a homogeneous porous medium. Assuming the absence of groundwater advection, the evolution of the soil temperature change $u$ in space $\mathbf{x}$ and time $t$ induced by a field of vertical BHEs is governed by the classical heat diffusion equation \cite{carslaw1948conduction}:
\begin{equation}
    \partial_t u(\mathbf{x}, t) - \thermaldif \Delta u(\mathbf{x}, t) = \frac{\thermaldif}{\thermalcond} f(\mathbf{x}, t)
\end{equation}
where $\thermaldif$ represents the thermal diffusivity and $\thermalcond$ denotes the thermal conductivity of the soil.

The source term $f$ models the boreholes as finite line sources of the length $\length$, positioned at coordinates $\positions$ in the horizontal plane. The heat injection/extraction rate per unit length, $q$, varies with time. To rigorously satisfy the semi-infinite geometry of the ground while enforcing a constant temperature Dirichlet boundary condition at the surface ($z=0$), we employ the Method of Images \cite{carslaw1948conduction}. This is represented in the source term by extending the domain to include a fictitious sink of equal magnitude but opposite polarity, located symmetrically above the ground surface ($z \in [-\length, 0]$):
\begin{equation}
    f(\mathbf{x}, t) = q(t) \left[ \mathbb{1}_{[0, \length]}(z) - \mathbb{1}_{[-\length, 0]}(z) \right] \delta_{\mathbf{p}}(x, y).
\end{equation}
\newpage
The fundamental analytical solution is derived by convolving this source distribution with the three-dimensional Green’s function (heat kernel) for infinite space \cite{evans2022partial}, given by:

$$\green(\mathbf{x}, t) = (4\pi\alpha t)^{-3/2} \exp\left(-\dfrac{|\mathbf{x}|^2}{4\alpha t}\right).$$  

By convolving the Green’s function with the vertical line source, we obtain a combination of error functions. We define the resulting vertical response function $\verticalkernel$ as:
\begin{equation} \label{def: vertical kernel}
    \verticalkernel(z, \tau;\length) := \frac{1}{2} \left[ \text{erf}\left(\frac{\length-z}{2\sqrt{\alpha\tau}}\right) + 2\text{erf}\left(\frac{z}{2\sqrt{\alpha\tau}}\right) - \text{erf}\left(\frac{\length+z}{2\sqrt{\alpha\tau}}\right) \right]
\end{equation} 
and the corresponding horizontal response function $\mathcal{R}$ as:
\begin{equation} \label{def: horizontal kernel}
    \mathcal{R}(r; \tau) :=\frac{1}{\tau} \exp\left( -\frac{r^2}{4\alpha\tau} \right).
\end{equation}

Consequently, the temperature variation at any position $\mathbf{x}=(x,y,z)$ and time $t$ is expressed as the temporal convolution of the load history with these spatial kernels:
\begin{equation}
    u(\mathbf{x}, t) = \frac{1}{4\pi\lambda} \int_{0}^{t} q(t-\tau) \cdot  \mathcal{R}(r(x,y, \positions); \tau) \mathcal{Z}(z, \tau; L) \, d\tau 
\end{equation}
where $\positions = (p_x, p_y), \quad r(x,y; \positions) = \| (x,y) - (p_x, p_y) \|$ is the radial distance from the borehole axis.

For the specific case where the heat load $q$ is constant or a step function, this temporal convolution is equivalent to the spatial integration form commonly found in literature, named as the Finite Line Source model \cite{LAMARCHE2007188, MOLINAGIRALDO20112506}:
\begin{equation}
u(\mathbf{x}, t) = \frac{q}{4 \pi  \lambda } \left(\int_0^L \frac{\text{erfc}\left(\frac{R(x,y,z; \xi)}
   {\sqrt{4 \alpha 
   t}}\right)}{R(x,y,z; \xi)} \,
   d\xi-\int_{-L}^0\frac{\text{erfc}\left(\frac{R(x,y,z; \xi)}
   {\sqrt{4 \alpha 
   t}}\right)}{R(x,y,z; \xi)} \,
   d\xi\right) 
\end{equation} 
with the Euclidean distance $ R(x,y,z; \xi) = \sqrt{r(x,y)^2 + (z-\xi)^2} .$ 
It is important to note that the integrands exhibit a singularity as the radial distance $r(x,y) \to 0$; therefore, the solution is valid strictly for evaluation points located at a non-zero radius, such as the borehole wall ($r_b > 0$), rather than at the centerline.

Due to the linearity of the governing heat equation w.r.t. the parameter $q$, the principle of superposition applies. This allows us to extend the single-source solution to a field of boreholes by summing the contributions of individual Dirac sources $\sum_{k=1}^{\nbhe} \delta_{\positions^k}$.
Strictly speaking, a coupled system would result in a distinct heat extraction rate $q_k$ for each borehole. However, we approximate the heat extraction rate as uniform across all boreholes ($q_k \approx q$), neglecting the slight variations that arise in practice. 
Hence, the total temperature field variation is obtained by summing over all borehole positions $\{\positions^k| k = 1, \dots N_b\}$:
\begin{equation} \label{eq: soil_temperature}
    u(\mathbf{x}, t) = \frac{1}{4\pi\lambda}\sum \limits _{k=1}^{\nbhe} \int_{0}^{t} q(t-\tau) \cdot  \mathcal{R}(r(x,y; \positions^k); \tau) \mathcal{Z}(z, \tau; L) \, d\tau .
\end{equation}

\subsection{Fluid Temperature Model}

We treat the depth range $z\in [0, L]$ as the active heat exchanger zone to model the boundary condition at the borehole wall. We use a quasi-steady-state approximation to describe how the circulating fluid and the borehole wall (at temperature $T_b(z,t)$) interact thermally. Neglecting the thermal mass of the fluid, the heat balance for the downward flowing fluid ($T_i$) and upward flowing fluid ($T_o$) is dictated by the coupled system of differential equations formulated by \cite{Eskilson01031988} for the 1U configuration:

\begin{equation}\label{eq: local_problem}
\begin{aligned}
    \density^{r}c^{r}\flowrate^r\frac{\partial T_{i}}{\partial z}(z,t) &= \frac{T_{b}(z,t) - T_{i}(z,t)}{R_{s}}-\frac{T_{i}(z,t)-T_{o}(z,t)}{R_{inter}}\\
   -\density^{r}c^{r}\flowrate^r\frac{\partial T_{o}(z,t)}{\partial z} &= \frac{T_{b}(z,t) - T_{o}(z,t)}{R_{s}}+\frac{T_{i}(z,t)-T_{o}(z,t)}{R_{inter}},
\end{aligned}
\end{equation}
with density $\density^r$, specific heat capacity $c^r$ and mass flow rate $\flowrate$. There are two main thermal resistances in the system $R_{s}$, which controls how heat moves between the pipe and the borehole wall, and $R_{inter}$, which accounts for how heat moves between the two pipe legs. These resistances can be determined experimentally or calculated using approximate analytical expressions based on the borehole geometry, as detailed in \cite{DIERSCH20111122, Diersch2014FEFLOW}.

\begin{figure}[h]
   \centering
   \includegraphics[width=0.5\textwidth]{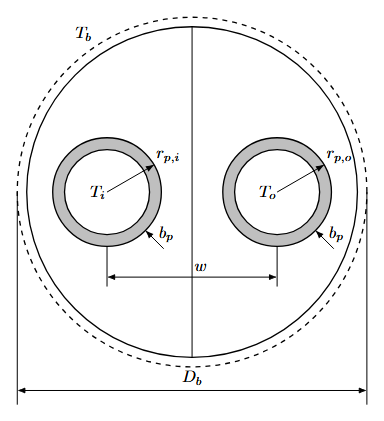}
   \caption{Cross section of a 1U borehole consisting of two pipe components and two grout zones (modified from \cite{DIERSCH20111122}).}
   \label{fig: 1U BHE}
\end{figure}

This formulation assumes thermal equilibrium within the borehole cross-section (Fig. \ref{fig: 1U BHE}). Thus, this description is only true for times $t>5\frac{\bheradius^2}{\thermaldif}$ \cite{DIERSCH20111122}, which gives the grout enough time to stabilize after any temporary changes.

Applying boundary conditions at the top and bottom of the borehole closes the equation system. We set a Dirichlet condition at the inlet ($z=0$) and a continuity condition at the U-bend ($z=L$) for a given time-dependent inlet temperature $T_{in}(t)$:
\[
\begin{aligned}
 & T_{i}(0,t)=T_{in}(t) \\
 & T_{i}(L,t)=T_{o}(L,t).
\end{aligned}
\]
It is assumed that the pipes are arranged symmetrically within the borehole for standard 1U configurations. This means that both legs have the same soil-to-pipe resistance $R_s$. The primary quantity of interest is the fluid outlet temperature $T_{out}(t) := T_{o}(0, t)$. This is the temperature of the fluid leaving the borehole at the surface.

The analytical solution to this system of equations \eqref{eq: local_problem}, which connects the temperature at the outlet to the temperature at the inlet and the vertical profile of the wall temperature, is found in \cite{Eskilson01031988}:
\begin{equation} \label{eq: outlet_temperature}
    T_{out}(t) = T_{in}(t)g_1(L) + g_2(L)\int_0^LT_b(\xi, t) g_3(L-\xi) d\xi
\end{equation}


The functions $g_1, g_2,$ and $g_3$ represent the BHE's geometric and thermal properties, which are defined as follows:
\begin{equation} \label{def: g-functions}
    \begin{aligned}
        g_1(z) &:= \frac{2 \gamma}{\beta_s \tanh (\gamma z)+\gamma}-1, \quad & 
        g_2(z) &:= \frac{\gamma}{\beta_s \sinh (\gamma z)+\gamma \cosh (\gamma z)}, \quad & 
        g_3(z) &:= 2 \beta_s \cosh (\gamma z) \\
        \text{where } \beta_s &= \frac{1}{R_{s}\density^{r}c^{r}\flowrate^r}, \quad & 
        \beta_{inter} &= \frac{1}{R_{inter}\density^{r}c^{r}\flowrate^r}, \quad & 
        \gamma &= \sqrt{\beta_s (\beta_s+2 \beta_{inter} )}
    \end{aligned}
\end{equation}


\subsection{Soil-Fluid Coupling strategy}

The system is driven by a given external energy demand $E$ varying in time. The coupling with the fluid loop is defined by the balance equation:
\begin{equation} \label{eq: balance}
    E(t) = \dot{m} c_p \left( \langle T_{out} \rangle(t) - T_{in}(t) \right)
\end{equation}
where $\langle T_{out} \rangle$ denotes the averaged outlet temperatures among the BHEs.
 We assume a uniform heat extraction rate along the length of each borehole and that these rates are identical for all boreholes in the field. But the average temperatures along the borehole walls are assumed to be unequal. This is due to the varying thermal interference experienced at different spatial positions within the field. This configuration corresponds to BC-I as defined in \cite{CIMMINO2014641}.

The Dirichlet boundary condition at the BHE walls implies that the soil temperature $u$ at the interface $r=r_b$ dictates the boundary temperature $T_b$. The thermal influence of a source borehole $k$ on a target borehole $\hat k$
 is given by:
\begin{equation} \label{def: u influence}
    u_{k\to \hat{ k}}(z,t) := \frac{1}{4\pi \lambda} \int_0^t q(t-\tau) \mathcal{R}(r_{k ,\hat{ k}}, \tau)\mathcal{Z}(z, L, \tau) \, d\tau,
\end{equation}
This formulation is motivated by \eqref{eq: soil_temperature}.
When considering the self-influence ($k=\hat{k}$), the temperature is computed at the borehole wall radius $r_b$. Therefore, the radial distance $r_{k ,\hat{ k}}$ is defined by:
\begin{equation}
    r_{k ,\hat{ k}} = \begin{cases}
        \| \mathbf{p}^k - \mathbf{p}^{\hat k} \| & \text{if } k\neq \hat k, \\
        r_b & \text{if } k=\hat k,
    \end{cases}
\end{equation} where $\mathbf{p}^k$ is the position of the BHE $k.$

To link the soil temperature field to the fluid outlet temperature, we introduce the vertically integrated thermal response, $\hat{u}_{k\to \hat k}$ between BHE $k$ and BHE $\hat k$ defined as
\begin{equation} \label{def: u hat}
    \hat{u}_{k\to \hat k}(t;\length) := \int_0^L u_{k\to \hat{k}}(z, t) \, g_3(\gamma(L-z)) \, dz.
\end{equation}
Physically, this quantity represents the effective temperature influence on the fluid, weighted by the BHE's internal heat transfer characteristics $g_3$ defined in \eqref{def: g-functions} along the length of the pipe.

Substituting \eqref{def: u influence} into \eqref{def: u hat} and applying Fubini's theorem to interchange the temporal and spatial integration yields:
\begin{equation}
    \hat{u}_{k\to \hat k}(t;\length) = \frac{1}{4\pi \lambda}\int_0^t q(t-\tau)\mathcal{R}(r_{k, \hat k}, \tau) \hat{\mathcal{Z}}(\tau;\length) \, d\tau,
\end{equation}
where the spatial integral $\hat{\mathcal{Z}}$ admits the following closed-form solution:
\begin{equation}
\begin{aligned}
    \hat{\mathcal{Z}}(\tau;\length) &:= \int_0^L \mathcal{Z}(z, \tau;\length) \cdot g_3(\gamma(L-z)) \, dz \\
    &= \frac{\beta_s}{2 \gamma} \mathcal{H}\left(\xi(\tau),\eta(\tau)\right).
\end{aligned}
\end{equation}
Here, the auxiliary function $\mathcal{H}$ is defined as:
\begin{equation}
   \mathcal{H}(\xi, \eta) = 4 \cosh^2(\xi\eta) \mathcal{E}(\xi, \eta) - \mathcal{E}(2\xi, \eta) - \big(1 + 2\cosh(2\xi\eta)\big) \mathcal{E}(0, \eta),
\end{equation}
with the function $\mathcal{E}$ given by
\begin{equation}
    \mathcal{E}(\xi, \eta) = e^{\eta^2} \left[ e^{2\xi\eta} \operatorname{erf}(\xi + \eta) - e^{-2\xi\eta} \operatorname{erf}(\xi - \eta) \right]
\end{equation}
with dimensionless parameters
$$ \xi(t) = \frac{L}{2\sqrt{\alpha t}}, \quad \eta(t) = \gamma\sqrt{\alpha t}. $$


By the symmetry of the radial distance matrix, we have $\hat{u}_{k\to \hat k} = \hat{u}_{\hat k\to k}$ for all $k, \hat k =1,\dots, \nbhe$. Furthermore, assuming identical borehole geometries, the self-response is uniform across the field, therefore $\hat{u}_{k\to k} = \hat{u}_{1\to 1}$ for all $k=1 , \dots, N_b$.

\newpage
Combining the soil temperature response \eqref{eq: soil_temperature} with the outlet temperature definition \eqref{eq: outlet_temperature}, and solving for the inlet temperature $T_{in}$ in the energy balance equation \eqref{eq: balance}, the average outlet temperature of the BHE field is derived as:
\begin{equation}
    \langle T_{out} \rangle(t) = \psi_1(L) \left( \Psi_{inter}(t;\length) + \Psi_{self}(t;\length) \right) + \psi_2(L) \Delta T_{pipe}(t).
\end{equation}
The coefficients $\psi_1$ and $\psi_2$, and the transient pipe temperature drop $\Delta T_{pipe}$ are defined as:
\begin{equation}
    \psi_1(z):= \frac{\gamma \text{csch}(\gamma z)}{2\beta_s}, \quad 
    \psi_2(z):= \frac{1}{2} \left(\frac{\gamma \coth (\gamma z)}{\beta_s}-1\right), \quad 
    \Delta T_{pipe}(t) := \frac{E(t)}{\dot{m} c_p}.
\end{equation}
Finally, the self-response and averaged interaction terms are given by:
\begin{align} \label{def: Psi}
    \Psi_{self}(t;\length) &:= \hat{u}_{1\to1}(t;\length) = \frac{1}{4\pi \lambda}\int_{0}^t q(t-\tau) \mathcal{R}(r_b, \tau)\hat{\mathcal{Z}}(\tau;\length) \, d\tau, \\
    \Psi_{inter}(t;\length) &:= \frac{2}{N_{bhe}} \sum_{k=2}^{N_{bhe}} \sum_{\hat k=1}^{k} \hat{u}_{k\to \hat k}(t;\length).
\end{align}

\subsection{Temporal Superposition and Aggregation} \label{subsec: Temporal superposition}

The heat extraction rate per unit length $q$ is related to the total building energy demand $E$ and the field configuration
\begin{equation}
q(t) = -\frac{E(t)}{\length \nbhe}.
\end{equation}
In practical applications, $E$ is modeled as a piecewise constant function defined on a uniform temporal grid $t_0, \dots, t_{N_t}$, typically representing hourly or monthly intervals
\begin{equation}
E(t) = \sum_{l =0}^{N_t-1}E_l \cdot \mathbb 1_{[t_l, t_{l+1})}(t).
\end{equation}
To accurately capture peak load behavior and ensure system limits are not violated, the simulation requires a fine temporal resolution (e.g., hourly time steps). However, a brute-force calculation of the full temperature field is computationally prohibitive. For a typical simulation horizon of 20 years with hourly steps ($N_t \approx 175,200$) and a field of 25 BHEs, a naive approach would require evaluating the convolution integral over $10^8$ time steps, scaling with $O(N_t^2 \cdot N_{bhe}^2)$.

To overcome this, we decompose the total temperature response into the local self-interaction term $\Psi_{self}$ and a global field-interaction term $\Psi_{inter}$, both defined in \eqref{def: Psi}.

Both terms rely on the horizontal response function $\mathcal{R}$, which exhibits a local maximum at $\tau_{\text{max}} = \frac{r^2}{4\alpha}$.

This duration corresponds to the moment of peak thermal impact at a radial distance $r > 0$. At the borehole radius $r = r_b (\approx 75 mm)$, the thermal response is almost instantaneous. For larger distances, such as the minimum BHE spacing $r \geq \Delta_{\min} (\approx 5m)$, the soil's thermal inertia acts as a low-pass filter: high-frequency fluctuations (e.g., hourly peaks) are significantly dampened and delayed by the time they reach neighboring boreholes. However, while individual pulses are preserved, the continuous injection or extraction of heat leads to a long-term thermal accumulation. During long operational durations, these energy flows accumulate, leading to a gradual shift of the temperature distribution that eventually influences the performance of the entire BHE field.



Based on this spectral separation, we employ a dual-grid approach.
For $\Psi_{inter}$ we aggregate the load $q$ onto a coarse grid (e.g., monthly averages). The convolution is performed on this reduced set of time steps, and the resulting smooth temperature curve is interpolated back to the hourly grid.
For $\Psi_{self}$ on the other hand, we perform the convolution on the original fine grid (hourly) to capture peak responses accurately. However, assuming uniform geometry of each BHE, this calculation is performed only once rather than for every pair of BHEs.

Finally, both $\Psi_{self}$ and $\Psi_{inter}$ can be cast as discrete convolutions. By utilizing the Fast Fourier Transform (FFT) and the Convolution Theorem, we reduce the complexity to $O(N_t \log N_t)$. For a 20-year hourly simulation, this reduces the operation count by several orders of magnitude, enabling the optimization of the BHE depth described in the following Section.

\section{Optimization Methodology}

\subsection{Algorithm Overview}

Our optimization strategy takes advantage of the fact that spatial density and required borehole length are linked in a physical way. To keep the total drilling depth and installation cost as low as possible, it is best to give each borehole as much thermal storage space as possible.

We do this with a two-step algorithm:
\begin{itemize}[leftmargin=1.5cm]
\item[Step 1] The algorithm puts a set number of BHEs inside the polygon that defines the property. The goal is to make the distance between any two boreholes and between any borehole and the property border as large as possible.
\item[Step 2] 
 With the BHE positions fixed, we perform a PDE-constrained optimization to determine the minimum uniform length $L$ along the BHE field. Using the analytical solution for the heat equation based on the FLS model, we identify the smallest $L$ required to satisfy the prescribed thermal demand while ensuring that the resulting fluid temperatures remain within predefined operational limits throughout the entire simulation period.
\end{itemize}

\subsection{Heuristic Positioning with Lloyd`s algorithm}

The spatial distribution objective described in Step 1 is formally framed as an optimization problem within the domain of Voronoi tessellations. Specifically, we seek to achieve a Centroidal Voronoi Tessellation (CVT) \cite{CVT} over the restricted, potentially non-convex domain $\Omega$ representing the property boundary. By characterizing the problem through its objective function \eqref{eq:objective} and associated constraints \eqref{eq:constraints1}–\eqref{eq:constraints3}, we ensure the optimal dispersion of BHEs. With a set number of BHEs $\nbhe$ that represent the desired number of Voronoi cells and a compact non-empty domain $\Omega \subset \mathbb{R}^2$, the goal is to find a set of $\nbhe$ generators (centroids) $\{z_i\}_{i=1}^{\nbhe} \subset \Omega$ and their corresponding Voronoi cells $\{V_i\}_{i=1}^{\nbhe}$ that minimize the following integral:

\begin{align}
\min_{\mathbf{z}, \mathbf{V}}
&\quad \mathcal{J}(\mathbf{z}, \mathbf{V})
:= \sum_{i=1}^{\nbhe} \int_{V_i} \|z_i-x\|^2\,dx
\label{eq:objective}\\
& \text{s.t.} \nonumber\\
& z_i \in V_i \quad \text{for all } i = 1, \dots, \nbhe \label{eq:constraints1}\\
& V_i \subset \Omega \quad \text{for all } i = 1, \dots, \nbhe \label{eq:constraints2}\\
& \bigcup_{i=1}^{\nbhe} V_i = \Omega, \quad
V_i^\circ \cap V_j^\circ = \emptyset \text{ for } i \ne j. \label{eq:constraints3}
\end{align}

The cells $V_i$ are thus defined as the regions in $\Omega$ closest to $z_i$:
    $$ V_i = \{ x \in \Omega \mid \|z_i - x\|^2 \le \|z_j - x\|^2 \text{ for all } j \ne i \} $$

For any given cell $V_i$, the point $z_i$ that minimizes $\int_{V_i} \|z_i-x\|^2\,dx$ is uniquely identified as the geometric centroid of $V_i$, given by $z_i = \fint_{V_i} x \, dx. $

Lloyd's algorithm \cite{Lloyd1982} is an iterative algorithm employed to approximate CVT. It consists of two alternating steps:
\begin{itemize}[leftmargin=1.5cm]
    \item[Step 1] Given the current set of generators $\{z_i\}$, the domain $\Omega$ is partitioned into Voronoi cells $V_i$. This step inherently satisfies the domain confinement constraint $V_i \subset \Omega$ by restricting the Voronoi regions to $\Omega$.
    \item[Step 2] Each generator $z_i$ is updated to the geometric centroid of its respective cell $V_i$.
\end{itemize}

In many real-world situations, the geometric domain $\Omega$ may not be convex. Using Lloyd's algorithm for CVT directly in these kinds of domains is problematic because the geometric centroid of a restricted Voronoi cell $V_i$ does not always lie within $V_i$, which results in the constraint \eqref{eq:constraints1} being violated. To address this issue, we implement an intuitive adjustment: if a computed centroid \(c_i\) is determined to be outside the domain \(\Omega\), we project \(c_i\) onto \(\Omega\).  There always exists a closest point, but the closest point may not be unique. Even though there is a chance that the operation will not be unique, it is still considered "well-defined" in a practical sense. Libraries for computational geometry (e.g. shapely \cite{gillies2013shapely}) are made to find the closest point. If a centroid $c_i$ is outside of $\Omega$, projecting it to the nearest point in $\Omega$ will bring it back into the feasible area, which is in line with the rule $z_i \in \Omega$.

To numerically approximate the CVT of a continuous domain $\Omega$, it is discretised into a sufficiently large quantity of uniformly distributed sample points \cite{CVT}.

Let $\mathcal{P} = \{p_1, \dots, p_M\}$ be a dense set of sample points within $\Omega$. The objective function $\mathcal J$ of the optimization problem \eqref{eq:objective} is then approximated as

$$ \min_{\mathbf{z}, \mathbf{S}}   \tilde{\mathcal{J}}(\mathbf{z}, \mathbf{S})=\sum_{i=1}^{\nbhe} \sum_{p_j \in S_i} \|z_i-p_j\|^2, $$

where $S_i \subset \mathcal{P}$ are discrete sets of sample points forming a partition of $\mathcal{P}$. The constraints are then applied to this discrete formulation, requiring $z_i \in \Omega$.

The modified Lloyd's algorithm, adapted for this discrete setting and incorporating the centroid projection, is described in Algorithm \ref{alg:modified_lloyds_cvt}.

\begin{algorithm}
\caption{Modified Lloyd's Algorithm for CVT in a Non-Convex Domain}
\label{alg:modified_lloyds_cvt}
\footnotesize
\begin{algorithmic}[1] 
\Require Non-convex domain $\Omega \subset \mathbb{R}^d$, number of generators $N$, maximum iterations $K$
, convergence tolerance $\varepsilon>0$.
\Ensure Set of converged generators $\{z_1, \dots, z_N\}$.

\State \textbf{Discretization of $\Omega$:} Generate a dense set of $M$ uniformly distributed sample points $\mathcal{P} = \{p_1, \dots, p_M\}$ such that $p_j \in \Omega$ for all $j=1,\dots,M$.
\State \textbf{Initialization of Generators:} Randomly select $N$ initial generators $\{z_1^{(0)}, \dots, z_N^{(0)}\}$ from $\mathcal{P}$.
\For {$t = 0, 1, \dots, K-1$}
    \State \textbf{Assignment Step:}
    \State Initialize empty discrete Voronoi cells $S_i^{(t)}$ for $i=1,\dots,N$.
    \For {each sample point $p_j \in \mathcal{P}$}
        \State Determine the closest generator $z_k^{(t)}$ to $p_j$ using Euclidean distance:
        \State $k^* = \operatorname{argmin}_{k \in \{1,\dots,N\}} \|z_k^{(t)} - p_j\|^2$
        \State Assign $p_j$ to the corresponding discrete Voronoi cell: $S_{k^*}^{(t)} \leftarrow S_{k^*}^{(t)} \cup \{p_j\}$
    \EndFor
    \State \textbf{Update  Step:}
    \For {$i = 1, \dots, N$}
        \If {$S_i^{(t)}$ is not empty}
            \State Calculate the geometric centroid (mean) $c_i^{(t+1)}$ of $S_i^{(t)}$:
            \State $c_i^{(t+1)} = \frac{1}{|S_i^{(t)}|} \sum_{p_j \in S_i^{(t)}} p_j$
            \State \textbf{Projection onto $\Omega$:}
            \If {$c_i^{(t+1)} \notin \Omega$}
                \State $z_i^{(t+1)} \in \operatorname{argmin}_{x \in \Omega} \|c_i^{(t+1)} - x\|^2$ \Comment{Project $c_i^{(t+1)}$ to any of the nearest points in $\Omega$}
            \Else
                \State $z_i^{(t+1)} = c_i^{(t+1)}$
            \EndIf
        \Else
            \State $z_i^{(t+1)} = z_i^{(t)}$ \Comment{Generator remains in place if its cell is empty}
        \EndIf
    \EndFor
    \State \textbf{Convergence Check:}
    \If {$\sum_{i=1}^{N} \|z_i^{(t+1)} - z_i^{(t)}\|^2 < \varepsilon$}
        \State \textbf{break}
    \EndIf
\EndFor
\State \Return $\{z_1^{(t+1)}, \dots, z_N^{(t+1)}\}$
\end{algorithmic}
\end{algorithm}

\begin{figure}[!htbp] 
    \centering
    \begin{subfigure}[b]{0.32\textwidth} 
        \includegraphics[width=\textwidth]{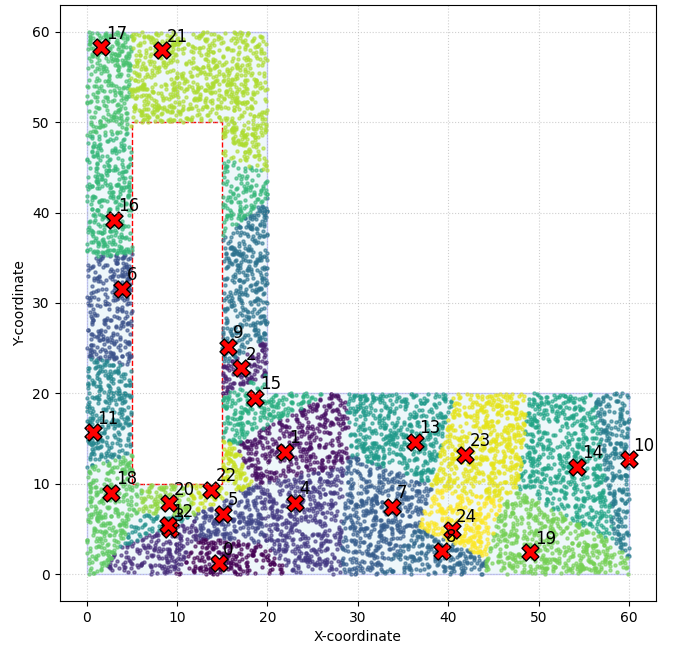} 
        \caption{Iteration 1}
        \label{fig:lloyds_iter1}
    \end{subfigure}
    \hfill 
    \begin{subfigure}[b]{0.32\textwidth}
        \includegraphics[width=\textwidth]{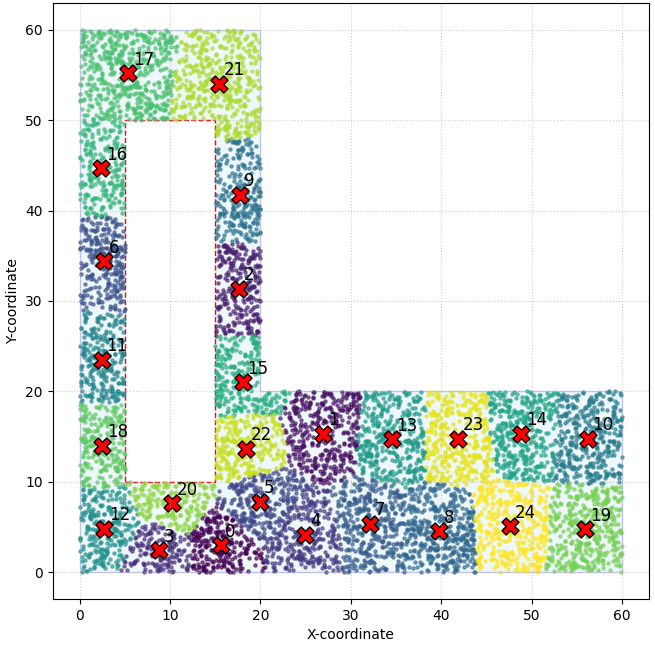} 
        \caption{Iteration 20}
        \label{fig:lloyds_iter20}
    \end{subfigure}
    \hfill
    \begin{subfigure}[b]{0.32\textwidth}
        \includegraphics[width=\textwidth]{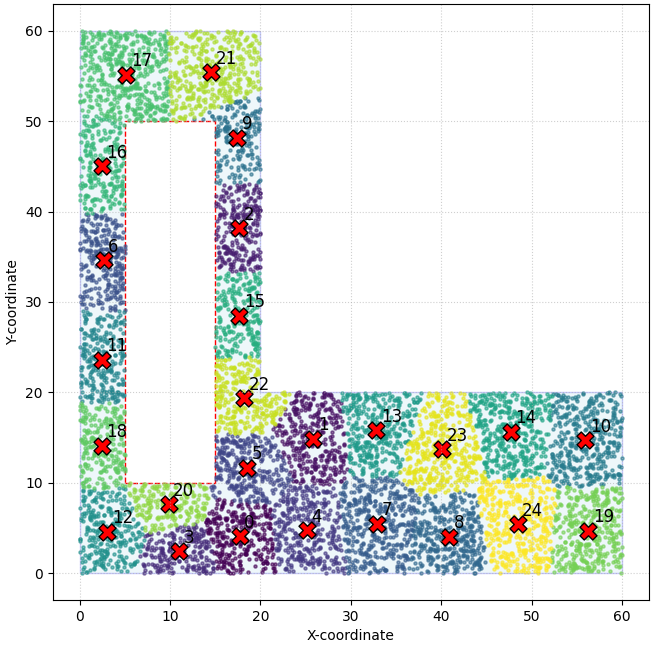} 
        \caption{Final State (Iteration 83)}
        \label{fig:lloyds_iterfinal}
    \end{subfigure}
    \caption{Progression of the Modified Lloyd's Algorithm in a non-convex domain $\Omega$. Sub-figure (a) shows the initial clustering, (b) demonstrates the state after 20 iterations as generators approach equilibrium, and (c) presents the final  converged configuration of the generators and their corresponding Voronoi cells.}
    \label{fig:lloyds_iterations}
\end{figure}

The iterative refinement process of the Modified Lloyd's Algorithm on the non-convex L-shaped domain $\Omega$ with a hole is illustrated in Figure \ref{fig:lloyds_iterations}, showcasing the evolution of generators and their associated Voronoi cells over several key steps.

\subsection{BHE Length Optimization}

The subsequent analysis focuses on optimizing the vertical dimension, building on the optimized spatial arrangement of BHEs that was generated applying the Modified Lloyd's Algorithm in the previous section. Specifically, we aim to determine the optimal uniform depth for these pre-positioned BHEs.

Let us introduce the optimization problem:
\begin{align}
       \min _{L} &\quad  L \label{eq:objective_L}\\
       \text{s.t.} &\quad T_{\min}\leq T_{\text{out}}(t;L)\leq T_{\max} \quad \forall t\in [0, \mathcal{T}]
       \label{constr: Tmin}\\
       &\quad L\in [L_{\min} , L_{\max}] \label{constr: control}
\end{align}
The objective function \eqref{eq:objective_L} is a scalar-valued linear function that aims to minimize the total cost associated with BHE installation. The uniform borehole length for all BHEs is represented by the decision variable $L$.

The problem is subject to several constraints. The control constraint \eqref{constr: control} sets practical limits on the length of the borehole $L$. The state constraints in \eqref{constr: Tmin} are crucial for making sure that the system operates properly by keeping the outlet fluid temperature from the BHEs, $T_{\text{out}}$, within a certain range of temperatures, $[T_{\min}, T_{\max}]$, for all times $t \in [0, \mathcal{T}]$. For numerical implementation, this continuous-time constraint is discretized over a finite set of time steps $\{t_l\}_{l=0}^{N_t-1}$ spanning the simulation period $\mathcal{T}$, leading to the practical form
$$
    T_{\text{out}}(t_l; L)\leq T_{\max} \quad \text{and} \quad T_{\text{out}}(t_l; L)\geq T_{\min} \quad \forall l=0, \dots, N_t-1.
$$
As explained in Section \ref{subsec: Temporal superposition}, the computation of $T_{\text{out}}$ implicitly takes into account the building's thermal demand profile. This demand determines the dynamic heat extraction or injection rates from the BHE field, which directly affects the temperature of the fluid flowing out of the BHE field. 

This optimization problem, featuring a linear objective, bound constraints, and non-linear state constraints, can be solved by various numerical optimization solvers. We use the Sequential Least Squares Programming (SLSQP) algorithm in Python's scipy.optimize \cite{Scipy} framework for our work. This choice is based on the fact that it performs well for medium-sized problems, handles both bound and non-linear constraints, and can be easily integrated to an open source development environment, making it a simple procedure to find a local optimum.

\section{Case Study}

To demonstrate our proposed heuristic approach, we present a case study with two scenarios. The thermal building loads for both scenarios are derived from a representative building model included in the TESS library \cite{TESS_Libraries}.
The remaining parameters characterizing the fluid, soil, and borehole geometry are adopted from the literature \cite{bernier2006closed, ahmadfard2019review}, and are comprehensively detailed in Table \ref{tab:fluid_soil_bhe_properties}. We will present the robustness of the heuristic approach by considering different property geometries, specifically the non-convex nature of the available land area for BHE placement.

\begin{table}[!tb]
    \scriptsize
    \centering
    \caption{Key Fluid, Soil, and Borehole Geometry Properties for BHE Simulation}
    \label{tab:fluid_soil_bhe_properties}
    \sisetup{
        separate-uncertainty = true,
        multi-part-units = single,
    }
    \begin{tabular}{llS[round-mode=places, round-precision=2]l} 
        \toprule
        Property & Symbol & {Value} & Unit \\
        \midrule
        \multicolumn{4}{l}{\textbf{Fluid Properties}} \\ 
        \midrule
        Specific Heat Capacity & $c_p$ & 4019 & \si{\joule\per\kilogram\per\kelvin} \\
        Density & $\rho$ & 1026 & \si{\kilogram\per\cubic\meter} \\
        Volumetric Heat Capacity & $c_v$ & \num{4.1245e6} & \si{\joule\per\cubic\meter\per\kelvin} \\
        Mass Flow Rate & $\dot{m}$ & \num{10.3397} & \si{\kilogram\per\second} \\
        Volume Flow Rate & $\dot{V}$ & \num{0.010078} & \si{\cubic\meter\per\second} \\
        Thermal Conductivity & $\lambda_f$ & 0.468 & \si{\watt\per\meter\per\kelvin} \\
        Dynamic Viscosity & $\mu$ & 0.00337 & \si{\kilogram\per\meter\per\second} \\
        \midrule
        \multicolumn{4}{l}{\textbf{Soil Properties}} \\ 
        \midrule
        Boundary Temperature & $T_s$ & 15 & \si{\degreeCelsius} \\
        Initial Inlet Temperature & $T_{in,0}$ & 15 & \si{\degreeCelsius} \\
        Thermal Conductivity & $\lambda_s$ & 1.9 & \si{\watt\per\meter\per\kelvin} \\
        Thermal Diffusivity & $\alpha_s$ & 0.08 & \si{\square\meter\per\day} \\
        Volumetric Heat Capacity & $C_s$ & \num{2.052e6} & \si{\joule\per\cubic\meter\per\kelvin} \\ 
        \midrule
        \multicolumn{4}{l}{\textbf{Borehole Geometry}} \\ 
        \midrule
        Borehole Diameter & $D_b$ & 150 & \si{\mm} \\ 
        Pipe Distance (center-to-center) & $w$ & 83 & \si{\mm} \\ 
        Outer Radius Pipe (in) & $r_{p,i}$ & 16.7 & \si{\mm} \\ 
        Outer Radius Pipe (out) & $r_{p,o}$ & 16.7 & \si{\mm} \\ 
        Pipe Wall Thickness & $b_{p}$ & 3.7 & \si{\mm} \\ 
        \bottomrule
    \end{tabular}
\end{table}

\begin{figure}[!htbp]
    \centering
    \includegraphics[width=0.8\textwidth]{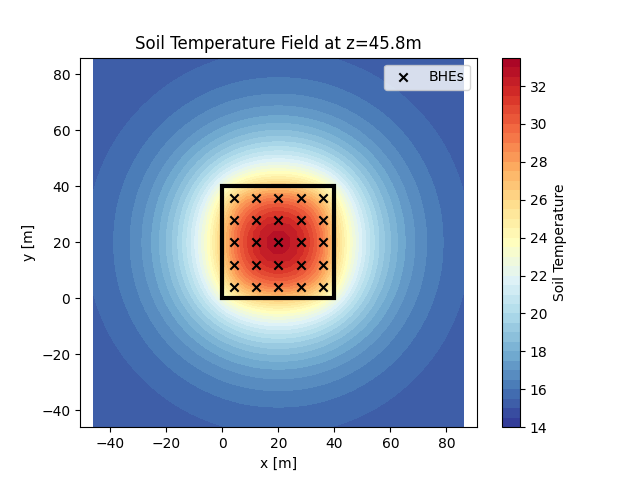}
    \caption{Soil temperature distribution at BHE mid-depth after a 20-year simulation period. Black crosses denote the positions of the BHEs, and the rectangle indicates the property boundary. The temperature field exhibits symmetric diffusion, consistent with the symmetry of the BHE array and the property borders, with the highest thermal accumulation concentrated around the centroid of the BHE field.}
    \label{fig:comparison_soil}
\end{figure}
\begin{figure}[htbp]
     \centering
     \begin{subfigure}[b]{0.32\textwidth}
         \centering
         \includegraphics[width=\linewidth]{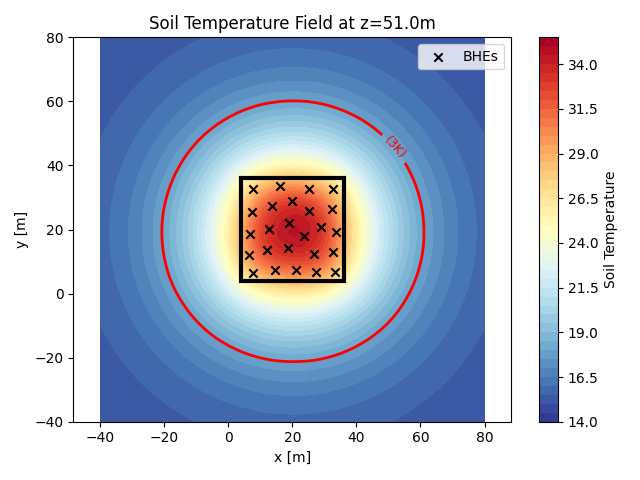}
         \caption{$1024\,\text{m}^2$ (Small)}
         \label{fig:first}
     \end{subfigure}
     \hfill
     \begin{subfigure}[b]{0.32\textwidth}
         \centering
         \includegraphics[width=\linewidth]{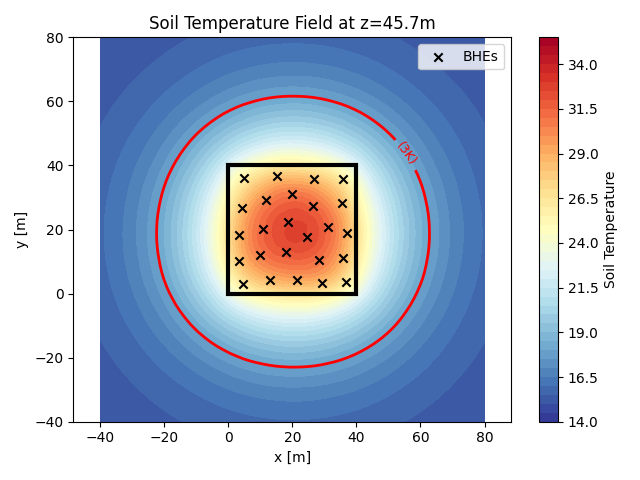}
         \caption{$1600\,\text{m}^2$ (Medium)}
         \label{fig:second}
     \end{subfigure}
     \hfill
     \begin{subfigure}[b]{0.32\textwidth}
         \centering
         \includegraphics[width=\linewidth]{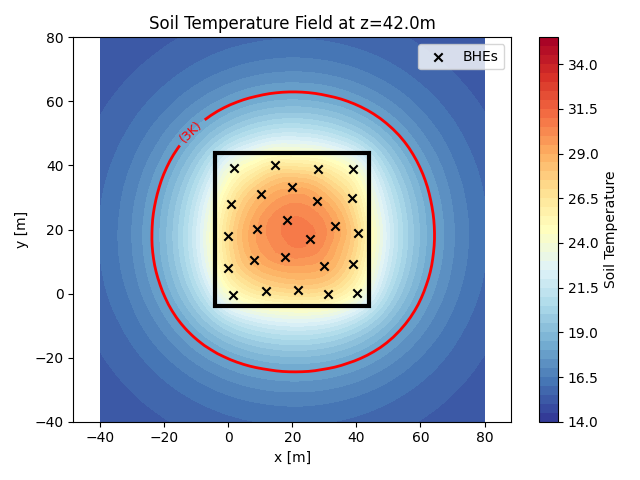}
         \caption{$2304\,\text{m}^2$ (Large)}
         \label{fig:third}
     \end{subfigure}
     \caption{Soil temperature distribution at BHE mid-depth for varying property dimensions. A clear correlation is observed between property size and thermal saturation, with peak temperatures increasing as the property area is constrained. In all cases, the red contour line delineates a temperature increase of $3\,\text{K}$ relative to the initial ground temperature.}
     \label{fig:three_graphs}
\end{figure}

We define a benchmark configuration that will serve as a baseline for comparison. The setup has 25 BHEs arranged in a standard 5 x 5 grid with 8 m of space between each of them on a 1600 $m^2$ rectangular domain. This symmetrical layout is typical for a manual setup of sizing tools. 
The performance of the heuristic approach is first evaluated by comparing the benchmark $5 \times 5$ grid (Fig. \ref{fig:comparison_soil}) with the positions generated by the Modified Lloyd's Algorithm for the equivalent $1600\,\text{m}^2$ domain (Fig. \ref{fig:second}) over a 20-year simulation period ($N_t = 175,200 = 8760\,\text{h/y} \times 20\,\text{y}$). The proximity between these configurations is significant, with an average distance of only $1.96\,\text{m}$ between the heuristically placed BHEs and their closest neighbors in the benchmark grid. This strong correlation indicates that under symmetric geometric constraints, the algorithm inherently reconstructs a quasi-grid structure, validating its ability to recover near-optimal placements in regular domains.

\begin{figure}[!htbp]
    \centering
    \includegraphics[width=0.75\linewidth]{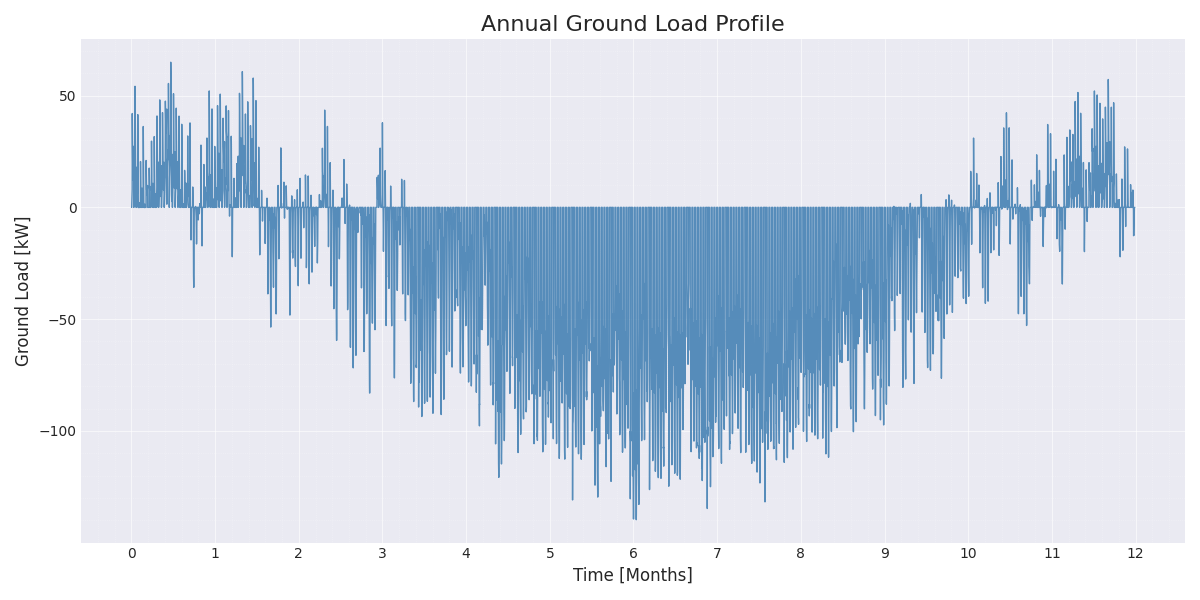}
    \caption{Annual Ground Load Profile with Hourly Discretization. This profile indicates how much thermal energy the BHE system must exchange at each hourly time interval over the period of a year. Positive values mean that heat is being extracted (heating demand), while negative values mean that heat is being rejected (cooling demand).}
    \label{fig: demand}
\end{figure}

\begin{figure}[!htb]
    \centering
    \includegraphics[width=0.75\linewidth]{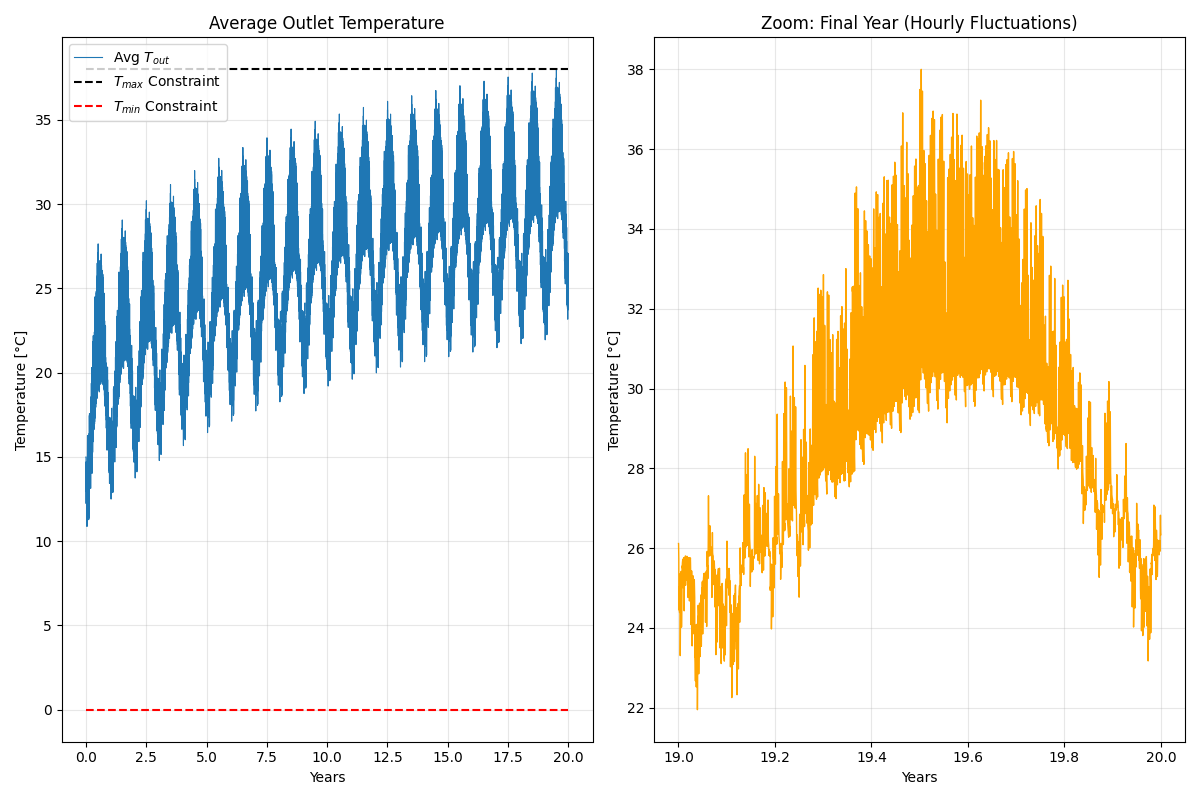}
    \caption{Evolution of the BHE outlet temperature over a 20-year horizon. The oscillating seasonal demand is superimposed on a rising mean temperature, eventually activating the $T_{\max}$ constraint.}
    \label{fig: Outlet temperature}
\end{figure}

The optimization results using SLSQP confirm this further: the optimal length for the manual grid is $L^*$ = $91.47\,\text{m}$, and the optimal length for the heuristic placement is $L^*$ = $91.44\,\text{m}$. Figure \ref{fig: demand} shows that this building's thermal demand is mostly for cooling. Because of this imbalance, the ground temperature rises over time, and the outlet temperature $T_{\text{out}}$ reaches the upper operational limit $T_{\max}$ at the end of the 20-year simulation period, as shown in Figure \ref{fig: Outlet temperature}.

The main benefit of the suggested heuristic lies in the fact that it can easily adjust to various domain constraints where manually fitting a grid is not straightforward. We examined the sensitivity of the optimal depth to the available property area by scaling the rectangle to $1024\,\text{m}^2$ (Fig. \ref{fig:first}) and $2304\,\text{m}^2$ (Fig. \ref{fig:third}).

The results, shown in Table \ref{tab:optimal_lengths_transposed_short}, indicate there is a clear inverse relationship between the amount of space available and the depth of the borehole required. The algorithm employs the extra space in the "Large" scenario to maximize BHE separation, thereby decreasing thermal interference and allows for a lower optimal depth of $84.07\,\text{m}$. On the other hand, the "Small" scenario needs a significant rise in depth to $101.86\,\text{m}$ to compensate for the higher thermal density. 

\begin{table}[!htbp]
    \scriptsize
    \centering
    \caption{Optimal BHE Length $L^*$ across various property geometries and scales (all lengths in \si{\meter}).}
    \label{tab:optimal_lengths_transposed_short}
    \sisetup{
        table-format=3.2,
        round-mode=places,
        round-precision=2
    }
    \begin{tabular}{l *{5}{S[table-format=3.2]}}
        \toprule
        & {Grid (5x5)} & {Rect. (Medium)} & {Rect. (Small)} & {Rect. (Large)} & {L-Shape} \\
        \midrule
        Optimal Length & 91.47 & 91.44 & 101.86 & 84.07 & 80.72 \\
        \bottomrule
    \end{tabular}
\end{table}

\begin{figure}[!htbp]
    \centering
    \includegraphics[width=0.75\linewidth]{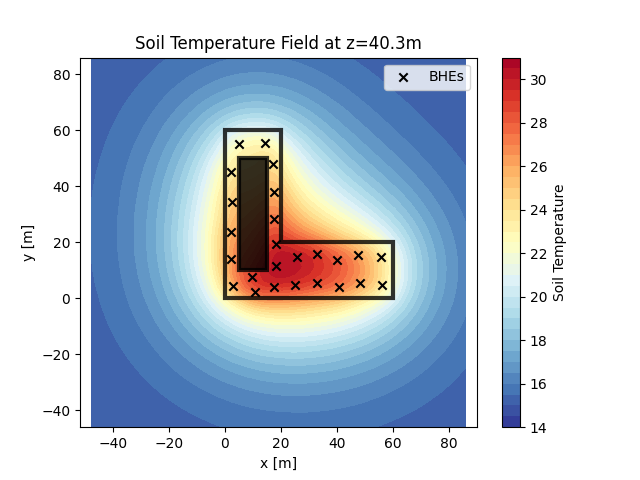}
    \caption{Soil temperature distribution in a non-convex L-shaped domain with an interior exclusion zone (building footprint). The algorithm automatically shifts BHE positions to optimize heat diffusion in the restricted space.}
    \label{fig: L-shaped soil temperature}
\end{figure}

To perform additional stress tests on the algorithm, we looked at a complicated L-shaped property with a rectangular exclusion zone that represents a building footprint. The Modified Lloyd's Algorithm effectively distributed the BHEs in the remaining non-convex area while considering the interior boundary. Figure \ref{fig: L-shaped soil temperature} displays the resulting soil temperature field. Despite the irregular geometry, the heuristic successfully identifies a configuration that maintains the fluid temperature within the required bounds, resulting in an optimal length of $80.72\,\text{m}$.

\section{Discussion and Outlook}

We developed a quick, coupled analytical framework which utilities a modified Lloyd's algorithm to find the optimal uniform depth for BHEs while also taking into account their spatial positioning. The heuristic effectively manages irregular and non-convex property boundaries, demonstrating its robustness in scenarios where traditional grid-based manual placement methods are impractical or inefficient. 
The main advantage of our method lies in its computational efficiency, enabling robust preliminary designs. The efficiency is achieved through a simple coupling of the Finite Line Source model with the inner 1D borehole heat transfer model, and the robustness through a modification of Lloyd´s algorithm for non-convex domains.
 Future research will explore gradient-based optimization methods for BHE coordinate positioning, going beyond the heuristic Lloyd's algorithm. This involves determining the gradients of the constraints in relation to BHE coordinates so that more rigorous gradient-descent algorithms can be used for spatial optimization. 
Another future improvement is to add explicit constraints on the maximum change in soil temperature ($\Delta T_{soil}$) to the optimization problem. It is crucial to keep the ground from freezing or overheating, which can lower system efficiency, damage the ground over time, and violate regulatory guidelines.

By including these feature in the futures, the presented framework can evolve into an even more powerful and accurate tool for the sustainable design of geothermal energy systems.

\section{Conflict of interest}
The authors declare no conflicts of interest. 

\section*{Acknowledgement}
\label{sec:Acknowledgement}
The work presented in this paper was supported by the German Federal Ministry of Economics and Climate Protection (BMWK) within the scope of the research project OptiGeoS (01256700/1).

\bibliography{bibliography.bib}

\end{document}